\pdfoutput=1
\RequirePackage{ifpdf}
\ifpdf 
\documentclass[pdftex]{sigma}
\else
\documentclass{sigma}
\fi

\begin{document}

\allowdisplaybreaks

\renewcommand{\thefootnote}{}

\renewcommand{\PaperNumber}{068}

\FirstPageHeading

\ShortArticleName{Numerical Approach to Painlev\'e Transcendents on Unbounded Domains}

\ArticleName{Numerical Approach to Painlev\'e Transcendents\\ on Unbounded Domains\footnote{This paper is a~contribution to the Special Issue on Painlev\'e Equations and Applications in Memory of Andrei Kapaev. The full collection is available at \href{https://www.emis.de/journals/SIGMA/Kapaev.html}{https://www.emis.de/journals/SIGMA/Kapaev.html}}}

\Author{Christian KLEIN and Nikola STOILOV}

\AuthorNameForHeading{C.~Klein and N.~Stoilov}
\Address{Institut de Math\'ematiques de Bourgogne, UMR 5584, Universit\'e de Bourgogne-Franche-Comt\'e,\\
 9 avenue Alain Savary, 21078 Dijon Cedex, France}
\Email{\href{mailto:Christian.Klein@u-bourgogne.fr}{Christian.Klein@u-bourgogne.fr}, \href{mailto:Nikola.Stoilov@u-bourgogne.fr}{Nikola.Stoilov@u-bourgogne.fr}}

\ArticleDates{Received April 18, 2018, in final form July 02, 2018; Published online July 12, 2018}

\Abstract{A multidomain spectral approach for Painlev\'e transcendents on unbounded domains is presented. This method is designed to study solutions determined uniquely by a, possibly divergent, asymptotic series valid near infinity in a sector and approximates the solution on straight lines lying entirely within said sector without the need of evaluating truncations of the series at any finite point. The accuracy of the method is illustrated for the example of the \emph{tritronqu\'ee} solution to the Painlev\'e~I equation.}


\Classification{34M55; 65L10}

\begin{flushright}
{\it This paper is dedicated to the memory of A.~Kapaev}
\end{flushright}

\renewcommand{\thefootnote}{\arabic{footnote}}
\setcounter{footnote}{0}

\section{Introduction}
Painlev\'e transcendents appear in many applications, see for instance \cite{Cla}, and have therefore been included as a chapter in the NIST digital library of mathematical functions~\cite{NIST}, thus complementing classical transcendents as hypergeometric functions. Consequently, the nume\-ri\-cal computation of these transcendents has seen a great deal of activity over the years. The goal of this paper is to offer an efficient numerical approach for transcendents characterized by their asymptotic behavior. The equations are studied on unbounded domains which offers a~boundaryless approach to Painlev\'e transcendents.

Numerical approaches to Painlev\'e transcendents found in the literature either use the integrability of the equations and essentially solve Riemann--Hilbert problems~\cite{FIKN}, see~\cite{olver} for numerical approaches, or evaluate Fredholm determinants~\cite{fredholm2,fredholm}, or they apply standard techniques for ordinary differential equations (ODEs). The latter allow also to study non-integrable deformations of the Painlev\'e equations with the same numerical techniques, e.g., to study cha\-rac\-teristic features of integrable ODEs not shared by non-integrable ones. These approaches can loosely be split into techniques solving initial value problems on one hand, see, e.g.,~\cite{JK}, and boundary value problems on the other, see, e.g., \cite{chebop,DubrovinGravaKlein,FW} and references therein. The boundaryless approach of the present paper is in particular intended to study the non-integrable deformations of Painlev\'e equations.

Many interesting transcendents are characterized by their asymptotic behaviour. The standard approach to numerically construct the transcendents is to provide a formal solution with the desired asymptotic in the form of an asymptotic series which is in general diverging. This series will be truncated according to standard rules, see, e.g.,~\cite{GR}. The truncated series will then provide boundary or initial data for the numerical solution at some finite value of the variable~$z$ in the complex $z$-plane. In this paper we address the obvious disadvantage of this approach, namely the use of diverging series and the arbitrariness in the choice of the value of~$z$ where to set up the initial or boundary value problem as well as in the truncation of the series. Instead, here we propose to study solutions determined by their asymptotics at infinity in some sector of the complex $z$-plane by numerically computing such solutions along a given straight line lying entirely within the sector with the use of the following techniques:
\begin{itemize}\itemsep=0pt
\item We use independent variable transformations to compactify the line into several finite intervals, such that the extreme endpoints of the collection of intervals correspond exactly to the points at infinity at the opposite ends of the straight line.

\item We regularise the desired solution at infinity by subtracting sufficiently many terms from its large $z$ asymptotic expansion so that the remainder tends to zero as $z \to \infty$ along the line. Then in each interval with an endpoint corresponding to $z\to \infty$ we solve for this normalised remainder rather than the Painlev\'e function itself. Of course, the remainder satisfies its own differential equation derivable explicitly from the Painlev\'e equation of interest.
\end{itemize}
The advantage of compactification is that we may use polynomial interpolation on finite intervals to numerically treat the boundary-value problem, and moreover working with the remainder gives the advantage that the boundary conditions are exactly known at the extreme endpoints of the collection of intervals. In this sense, our method is ``boundary-less''. Note that the equations studied near the infinite $z$-points are singular, and thus no boundary conditions will be applied there.

The example for which we illustrate the concept is the \emph{tritronqu\'ee} solutions to the Painlev\'e~I (PI) equation,
\begin{gather}
 \frac{d^{2}\Omega}{dz^{2}}=3\Omega^{2}-z, \qquad z\in \mathbb{C}, \label{PI}
\end{gather}
which are characterized by Boutroux \cite{BoutrouxI,Boutroux} by their asymptotic
behaviour,
\begin{gather}
 \Omega \sim \sigma \sqrt{\frac{z}{3}}, \qquad z\to\infty,\qquad \sigma=\pm1. \label{asym}
\end{gather}
Boutroux showed that there are solutions in four consecutive sectors of the complex plane of angle $2\pi/5$ with the behaviour at infinity as in~(\ref{asym}), which asymptotically are pole-free and which he called \emph{tritronqu\'ee}. The sector can be chosen to be given by $|\operatorname{arg}z|<4\pi/5$ and we use the notation $\Omega(z)$ in the rest of the paper to refer to this well-defined particular solution of~(\ref{PI}). There are four other \emph{tritronqu\'ee} solutions related to this solution via
\begin{gather*}
 \Omega_{n}(z)=\mathrm{e}^{4\pi\mathrm{i}n/5}\Omega\big(\mathrm{e}^{2\pi\mathrm{i}n/5}z\big) ,\qquad n=\pm1,\pm2, 
\end{gather*}
which corresponds to a general symmetry of the PI equation not limited to \emph{tritronqu\'ee} solutions. Note that the considered \emph{tritronqu\'ee} solution has the symmetry $\bar{\Omega}(z)=\Omega(\bar{z})$. It is shown in~\cite{ince} that any solution to the PI equation is a meromorphic function on the complex plane. For a~comprehensive discussion of \emph{tritronqu\'ee} solutions see Joshi and Kitaev~\cite{JK}. Kapaev~\cite{K1,K2,kap_p1_quasi} gave a complete characterization of the \emph{tritronqu\'ee} solutions in terms of the Riemann--Hilbert problem associated to the PI equation.

In~\cite{DubrovinGravaKlein} these solutions are applied to describe asymptotically the behaviour of solutions to the nonlinear Schr\"odinger (NLS) equation in the semiclassical limit near a critical point. It was conjectured that the \emph{tritronqu\'ee} solutions are pole-free in the whole sector $|\operatorname{arg}z|<4\pi/5$, not just asymptotically, which was then proven in~\cite{CHT}. In this paper we are interested in the \emph{tritronqu\'ee} solutions in the regular sector and will construct them on unbounded lines via polynomial interpolation. If poles are to be studied, approaches based on Pad\'e approximants as in \cite{FFW,FW, Nov} are to be used.

The paper is organized as follows: in Section~\ref{section2} we present the numerical techniques we employ. In Section~\ref{section3} we apply a multi-domain spectral approach to the \emph{tritronqu\'ee} solution on unbounded
lines. We add some concluding remarks in Section~\ref{section4}.

\section{Multidomain spectral method}\label{section2}

In this section we briefly summarize the numerical methods we use. We are interested in solutions on lines $z=ax+b$ in the complex plane, where $a,b\in\mathbb{C}$ are constant with respect to $x\in\mathbb{R}$. The line will be divided into several intervals for $x$ each of which will be mapped to the inter\-val~$[-1,1]$. On the latter, polynomial interpolation is used to approximate the solution.

As discussed for instance in \cite{FW}, it is numerically problematic to set up an initial value problem for the \emph{tritronqu\'ee} solution for large $|z|$, say for $x_{0}\gg1$ and to integrate the solution up to the value $-x_{0}$ since cancellation errors will lead to a loss of digits which can destabilize such a \emph{shooting} approach. In practice it is numerically more stable to set up a boundary value problem with boundary data at~$\pm x_{0}$.

\subsection{Domains and equations}
The line $z=ax+b$, $x\in\mathbb{R}$ will be divided into three domains, $x<x_{l}$, $x_{l}<x<x_{r}$ and $x>x_{r}$ numbered I, II, III respectively; here $x_{l}<x_{r}$ are real constants.

The \emph{tritronqu\'ee} solution we are interested in is characterized by the asymptotic beha\-viour~(\ref{asym}). This implies, see~\cite{JK}, that near infinity the solution can be represented by a formal series,
\begin{gather}
 \Omega(z)\sim\sigma \sqrt{\frac{z}{3}}+\sum_{k=1}^{\infty}\frac{a_{k}}{z^{(5k-1)/2}}, \quad\quad |\arg z| < 4\pi/5. \label{asymseries}
\end{gather}
The series in (\ref{asymseries}) is divergent and therefore of little use from an analytical point of view. However, its form suggests a numerical approach to solve for the function
\begin{gather}
 v:=\Omega-\sigma\sqrt{\frac{z}{3}},\qquad s := \frac{1}{\sqrt{ax+b}}. \label{vs}
\end{gather}
This function is regular at infinity in the sector under consideration: $v \to 0$ as $z \to \infty$ for $|\arg z| < 4\pi/5$. The PI equation (\ref{PI}) implies the following equation for $v$,
\begin{gather}
 \frac{s^{7}}{4}v_{ss}+\frac{3}{4}s^{6}v_{s}-2\sigma\sqrt{3}v= 3sv^{2}+\frac{\sigma}{4\sqrt{3}}s^{4}. \label{veq}
\end{gather}
Note that the same approach can be used for similar Painlev\'e transcendents though exponentially small terms in the asymptotics might require the use of more general transformations than considered here. Furthermore, starting with a function that is singular at infinity there are other ways to generate a function that vanishes there, e.g., in the considered case one could take~$1/\Omega(z)$ that would have a prescribed algebraic decay in $s$ at infinity. This approach will not be explored here.

This equation will be solved in domains I and III, equation (\ref{PI}) in domain II. Each of these domains will be mapped to the interval $[-1,1]$ in the following way ($l\in[-1,1]$),
\begin{gather*}
 s = \frac{1}{\sqrt{ax_{l}+b}}\frac{1+l}{2} 
\end{gather*}
in domain I,
\begin{gather*}
 x = x_{l}\frac{1-l}{2}+x_{r}\frac{1+l}{2} 
\end{gather*}
in domain II and
\begin{gather*}
 s = \frac{1}{\sqrt{ax_{r}+b}}\frac{1+l}{2} 
\end{gather*}
in domain III. Note that the choice of $x_{l}$ and $x_{r}$ can be optimized according to \cite{DW} (we mainly check this via the Chebyshev coefficients as shown in the next section), and that it is straight forward to generalize the approach to more than 3 domains. This was just chosen here to keep the presentation simple.

At the domain boundaries, the function $\Omega$ has to be $C^{1}$ in $x$. Together with the PI equation, this guarantees that the solution will be smooth on the whole considered line $x\in \mathbb{R}$. Since equation (\ref{veq}) is singular for $s=0$, i.e., for $z\to\infty$, no condition needs to be given there (the numerical approach will automatically produce the regular solution at this point).
\begin{remark} All roots appearing in this paper are to be understood as being defined on a two-sheeted Riemann surface. The sign of the roots will be fixed as usual at some point which is not a branch point. By analytic continuation as for instance in~\cite{FK}, sheets on the Riemann surface can be numerically established independently of the definition of the root in Matlab.
\end{remark}

\subsection[Polynomial interpolation and $\tau$-method]{Polynomial interpolation and $\boldsymbol{\tau}$-method}
\looseness=1 Since each of the above three domains have been mapped to the interval $[-1,1]$, the task is now to approximate a smooth function on this interval. A~standard approach is to choose a~suitable discretisation of the independent variable $l$, see for instance \cite{chebop, trefethen,WR}, and to use polynomial interpolation on the \emph{collocation points} $l_{j}$, $j=0,\ldots,N_{c}$. This leads to an approximation of the derivatives in terms of so-called differentiation matrices obtained by differentiating the interpolation polynomials. In~\cite{GK1}, a~collocation method with cubic splines distributed as \emph{bvp4} with Matlab was applied. In \cite{GK12}, a Chebyshev collocation method on Chebyshev collocation points $l_{j}=\cos(j\pi/N_{c})$, $j=0,\ldots,N_{c}$ was used. The latter is related to an expansion of the solution in terms of Chebyshev polynomials $T_{n}(l)=\cos(n\arccos(l))$, $n=0,1,\ldots$ of some func\-tion~$u(l)$,
\begin{gather}
 u(l)\approx\sum_{n=0}^{N_{c}}c_{n}T_{n}(l). \label{cn}
\end{gather}
Putting $y=\arccos(l)$, $a_{0}=c_{0}$, and $a_{n}=a_{-n}=c_{n}$ for $n=1,2,\ldots,N_{c}$, one gets for rela\-tion~(\ref{cn})
\begin{gather}
 u(l)\approx\sum_{n=-N_{c}}^{N_{c}}a_{n}e^{in y} \label{an},
\end{gather}
i.e., a discrete Fourier transform in the variable $y$ with collocation points $y_{j}=\pi j/N_{c}$, $j=-N_{c},\ldots,N_{c}$. Thus one advantage of the use of Chebyshev polynomials is that the expansion in terms of the latter can be efficiently computed with a \emph{fast cosine transform} (FCT) which is related via (\ref{an}) to a \emph{fast Fourier transform} (FFT) for which efficient algorithms exist (note that special care has to be taken of the terms $j=\pm N_{c}$ in order to establish the relation to an FFT, see the discussion in~\cite{trefethen}).

\looseness=-1 A further advantage of the relation between Chebychev and Fourier series is that results for the decrease of the Fourier coefficients of a given function $f$ for large $n$ in dependence of the regularity of $f$ can be applied directly to Chebychev series, see~\cite{trefethen,trefethen2}. Since Fourier~(\ref{cn}) and Chebychev sums~(\ref{an}) can be seen as truncated series, the numerical error in approximating a~given function by a~sum is of the order of the coefficients $c_{N_{c}}$ and $a_{N_{c}}$ respectively. We recall from~\cite{trefethen,trefethen2} that for a~function $f\in C^{p-1}(\mathbb{T})$ with a $p$th derivative of bounded variation, the Fourier coefficient $a_{n}=O\big(n^{-(p+1)}\big)$ for $n\to\infty$. For a~function $f\in C^{\infty}(\mathbb{T})$, the Fourier coefficients decrease as $a_{n}=O(n^{-m})$, $\forall\, m\in \mathbb{N}$. In numerical applications where one works with finite precision, this means that $C^{\infty}$ functions are approximated with an essentially exponentially small error. Note that because of the relation between Chebychev sums and discrete Fourier transforms, functions $u\in C^{\infty}([-1,1])$ are approximated with essentially exponential accuracy by the sum~(\ref{cn}). This allows to control the numerical resolution via the spectral coefficients.

\begin{remark}The above description of properties of Chebychev sums motivates our choice of the independent variable $s = 1/ \sqrt{z}$ for the remainder function $v$ defined in intervals I and III when we compute the Painlev\'e-I \emph{tritronqu\'ee} solution $\Omega(z)$. Indeed, the asymptotic expansion~(\ref{asymseries}) implies that $v(s)$ has an asymptotic expansion as $s \to 0$ with $| \arg(s)| < 2\pi/5$, and this expansion is a~(divergent) power series in s. Since $\Omega(z)$ is known to be analytic in the sector $| \arg(z)| < 4\pi/5$, it follows that $v(s)$ is analytic for $|s| > 0$ in the sector $| \arg(s)| < 2\pi/5$, and then the Cauchy integral formula shows that the asymptotic power series for $v(s)$ is differentiable term-by-term, i.e., arbitrary derivatives of $v(s)$ have asymptotic power series expansions as $s \to 0$ in the indicated sector given by
formal derivatives of the power series for $v$. In particular, this means that if the line $z = ax + b$ is taken with $| \arg(a)| < 4\pi/5$, all derivatives of $v(s)$ will be continuous up to the endpoints when intervals~I and~III are mapped to $[-1, 1]$, which guarantees that the Chebychev coefficients will satisfy $c_{n} = \mathcal{O}(n^ {-m} )$, $\forall\, m \in N$ as $n \to \infty$. However if $|\arg(a)|$ equals or exceeds $4\pi/5$ we can draw no such conclusion from the asymptotic series~(\ref{asymseries}).
\end{remark}

The PI equation (\ref{PI}) is thus replaced by $N_{\rm I}+1$, $N_{\rm II}+1$ and $N_{\rm III}+1$ algebraic equations where $N_{\rm I}$, $N_{\rm II}$, $N_{\rm III}$ give the number of collocation points on the respective domain. The resulting system of algebraic equations is solved using Newton's method. As the initial iterate we use in domains~I and~III $v=0$; for domain~II, we compute the first 6 terms of the asymptotic series~(\ref{asymseries}) for $x_{l}$ and $x_{r}$ and use the linear interpolate between these values as an initial iterate.

The normal Newton iteration for the solution of an equation $F(u)=0$ takes the form
\begin{equation*}
u_{n+1}=u_{n}-(\operatorname{Jac}F(u_{n}))^{-1}F(u_{n}),
\end{equation*}
where $\operatorname{Jac} F$ is the Jacobian of $F$ and $u_{n}$ the $n$th iterate. The accessible precision is mainly limited by the conditioning of the Chebyshev differentiation matrices which is of the order of~$N_{c}^{2}$, see the discussion in~\cite{trefethen} for second order differentiation matrices. Note that this problem can be addressed as mentioned by introducing more than 3 domains. The iteration is stopped typically at a residual $||F(u_{n})||_{\infty}$ of~$10^{-10}$.

The junction conditions at the domain boundaries are implemented via Lanczos' $\tau$-me\-thod~\cite{tau}. This means that the $C^{1}$ conditions are simply replacing the lines corresponding to $x=x_{l}$ and $x=x_{r}$ in the equation
\begin{gather*}
 \operatorname{Jac}F((u_{n}))(u_{n+1}-u_{n})=-F(u_{n}),
\end{gather*}
before numerically solving for $u_{n+1}$. The derivatives in the differentiability conditions on $u$ at the boundary are again computed via Chebyshev differentiation matrices. It is known, see~\cite{trefethen}, that the boundary conditions implemented in this way will be satisfied with the same spectral accuracy as the solution.

\section{Tritronqu\'ee solutions on lines in the complex plane}\label{section3}
In this section we consider two examples for the above numerical scheme, the \emph{tritronqu\'ee} solution to PI on the imaginary axis and a line close to the \emph{Stokes' lines} $\operatorname{arg}z=\pm 4\pi/5$ where the solution is known to show an oscillatory behavior, see Kapaev~\cite{kap_p1_quasi}. On the Stokes' lines, the \emph{tritronqu\'ee} solutions have an oscillatory singularity at infinity, and the spectral approach will obviously have problems to approximate such a behavior.

\subsection{Imaginary axis}
\begin{figure}[t!]\centering
 \includegraphics[width=0.49\textwidth]{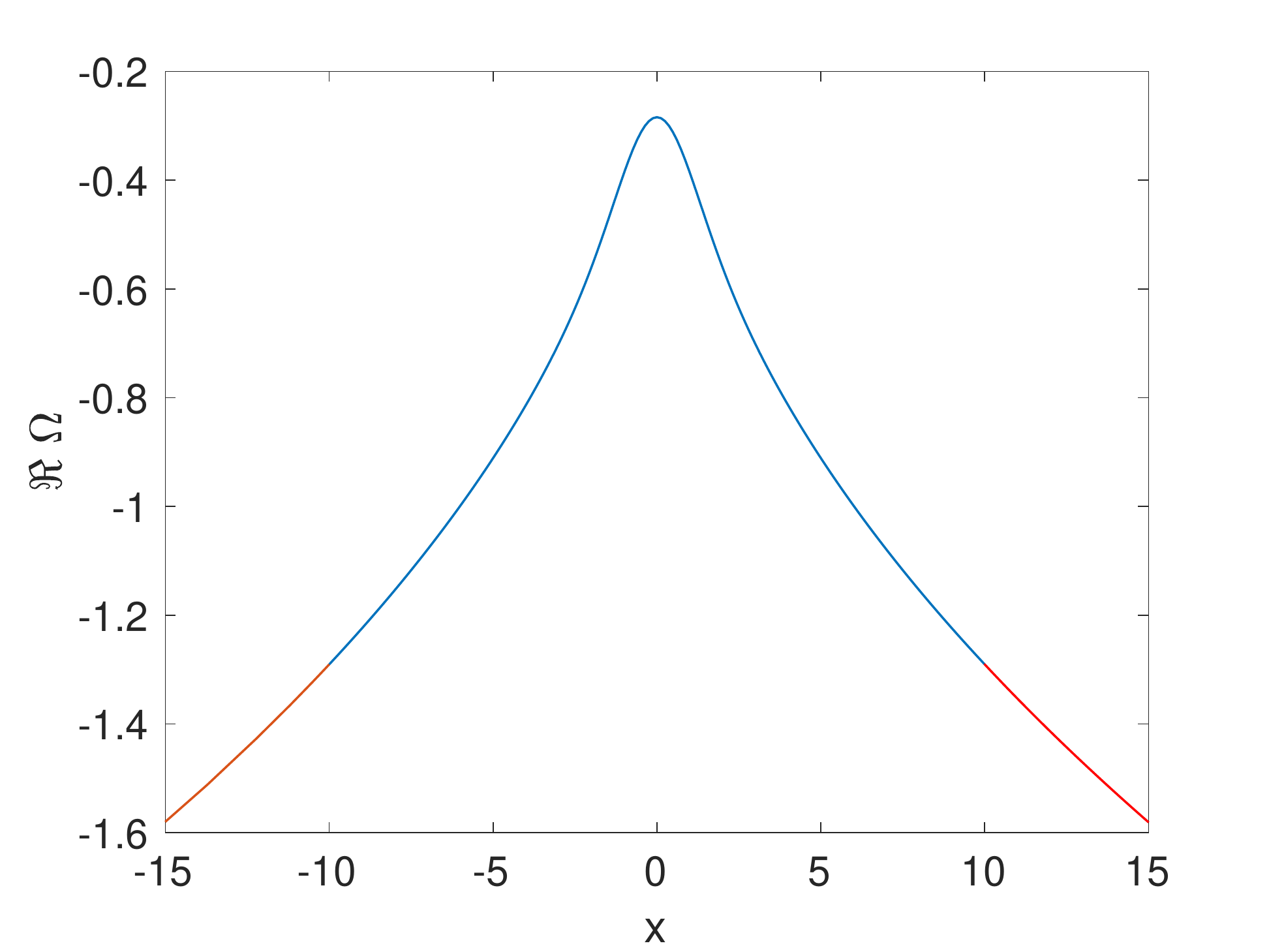}
 \includegraphics[width=0.49\textwidth]{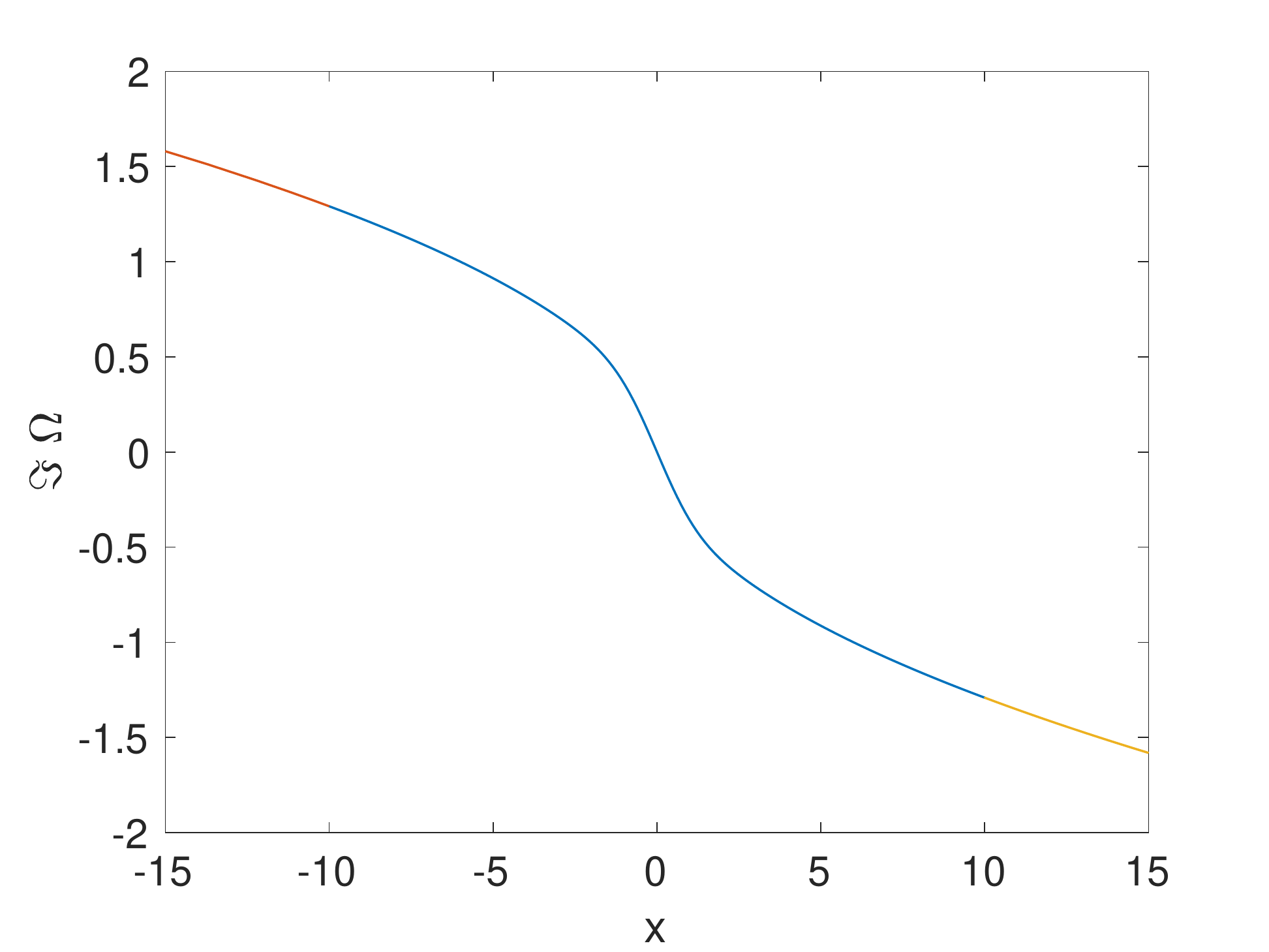}
 \caption{\emph{Tritronqu\'ee} solution to the PI equation on the imaginary axis $z=\mathrm{i}x$, on the left the real part, on the right the imaginary part; the solutions are shown in blue in domain~ II and in red in domains~I and~III.} \label{imag}
\end{figure}

As a first example, we will study the case $z=\mathrm{i}x$. We use $x_{r}=-x_{l}=10$ and $N_{\rm I}=N_{\rm III}=20$ and $N_{\rm II}=256$ collocations points in the respective domains. The solution $\Omega$ is shown for $|x|<15$ in Fig.~\ref{imag}. It can be recognized that the solution is smooth on the whole interval and in particular at the domain boundaries $x_{l}$, $x_{r}$. The symmetry $\bar{\Omega}(z)=\Omega(\bar{z})$ of the solution is obvious.

\begin{figure}[t!]\centering
 \includegraphics[width=0.49\textwidth]{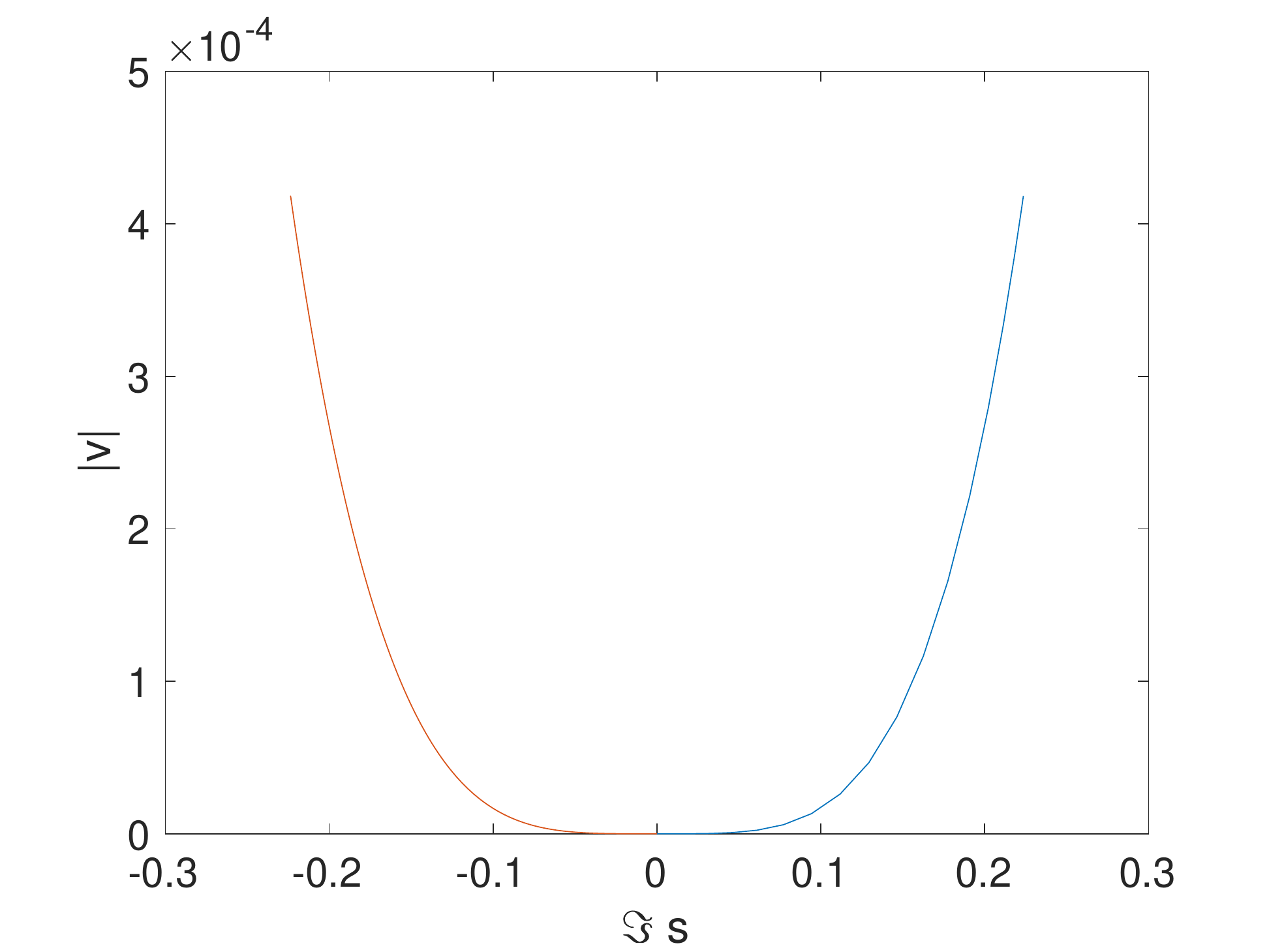}
 \includegraphics[width=0.49\textwidth]{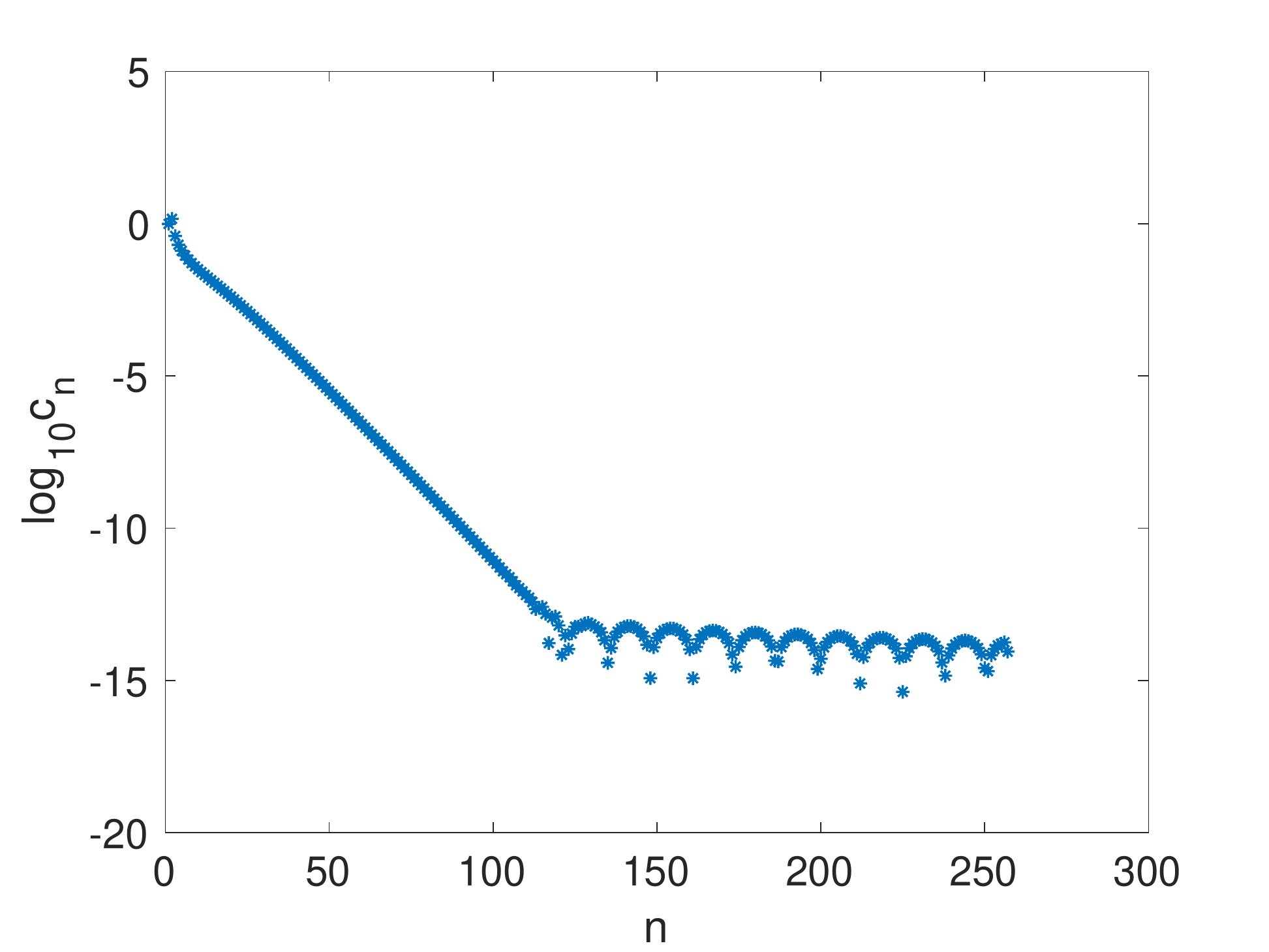}
 \caption{On the left, the computed functions $v$ corresponding to Fig.~\ref{imag} in domains I (blue) and III (red); on the right the modulus of the spectral coefficients $c_{n}$ in domain II in a~logarithmic plot.}
 \label{sfig}
\end{figure}

The function $\Omega$ in domains I and III follows from the computed functions $v$ via (\ref{vs}). The computed $v$ in dependence of~$s$, i.e., on the unbounded domains, is shown on the left of Fig.~\ref{sfig}. It can be seen that the functions are of the order of $10^{-4}$ and smooth on the considered intervals. As expected they vanish at infinity. This shows that the \emph{tritronqu\'ee} solution is in the considered example already very close to the asymptotic solution.

\begin{figure}[!t]\centering
 \includegraphics[width=0.49\textwidth]{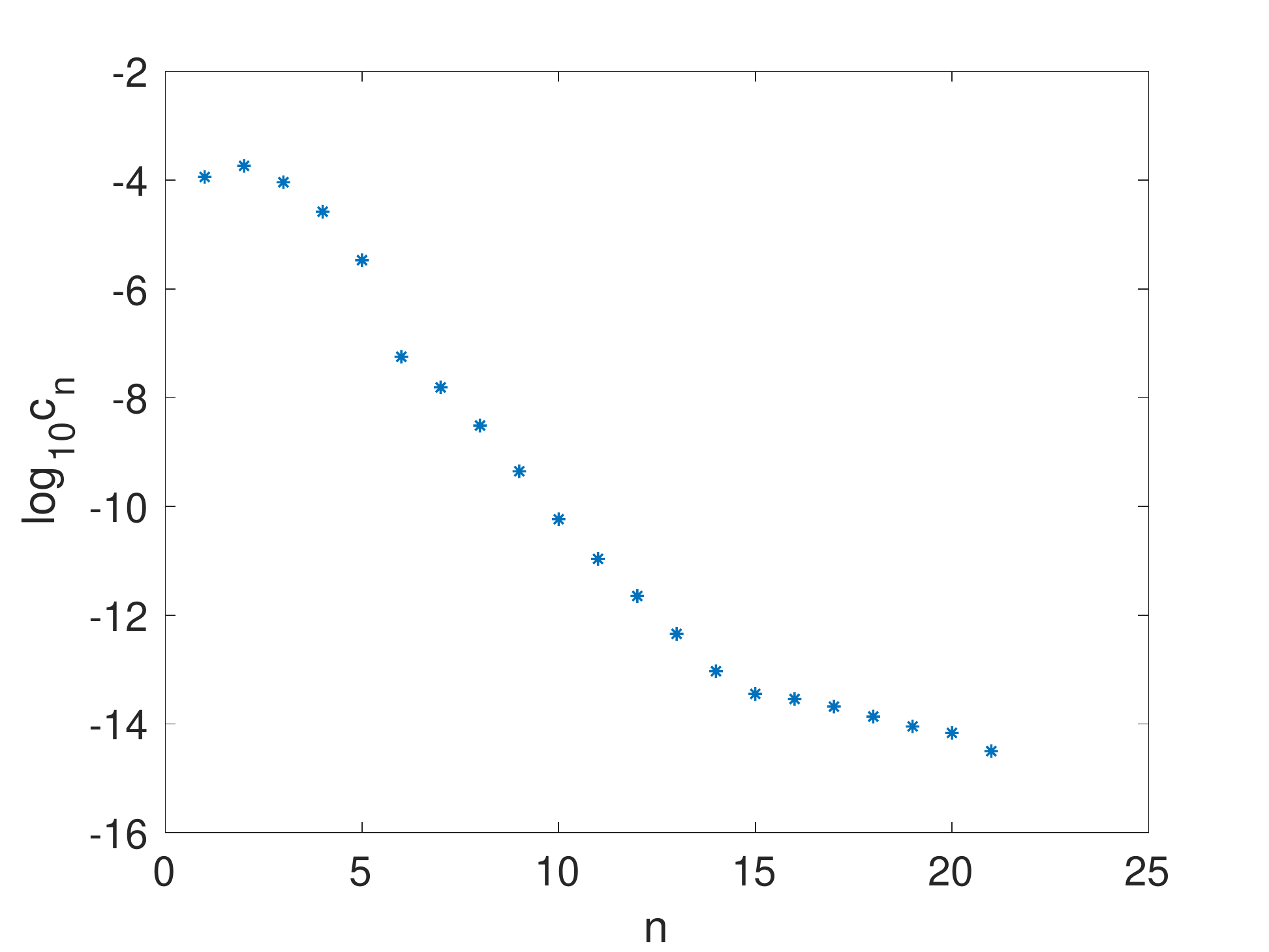}
 \includegraphics[width=0.49\textwidth]{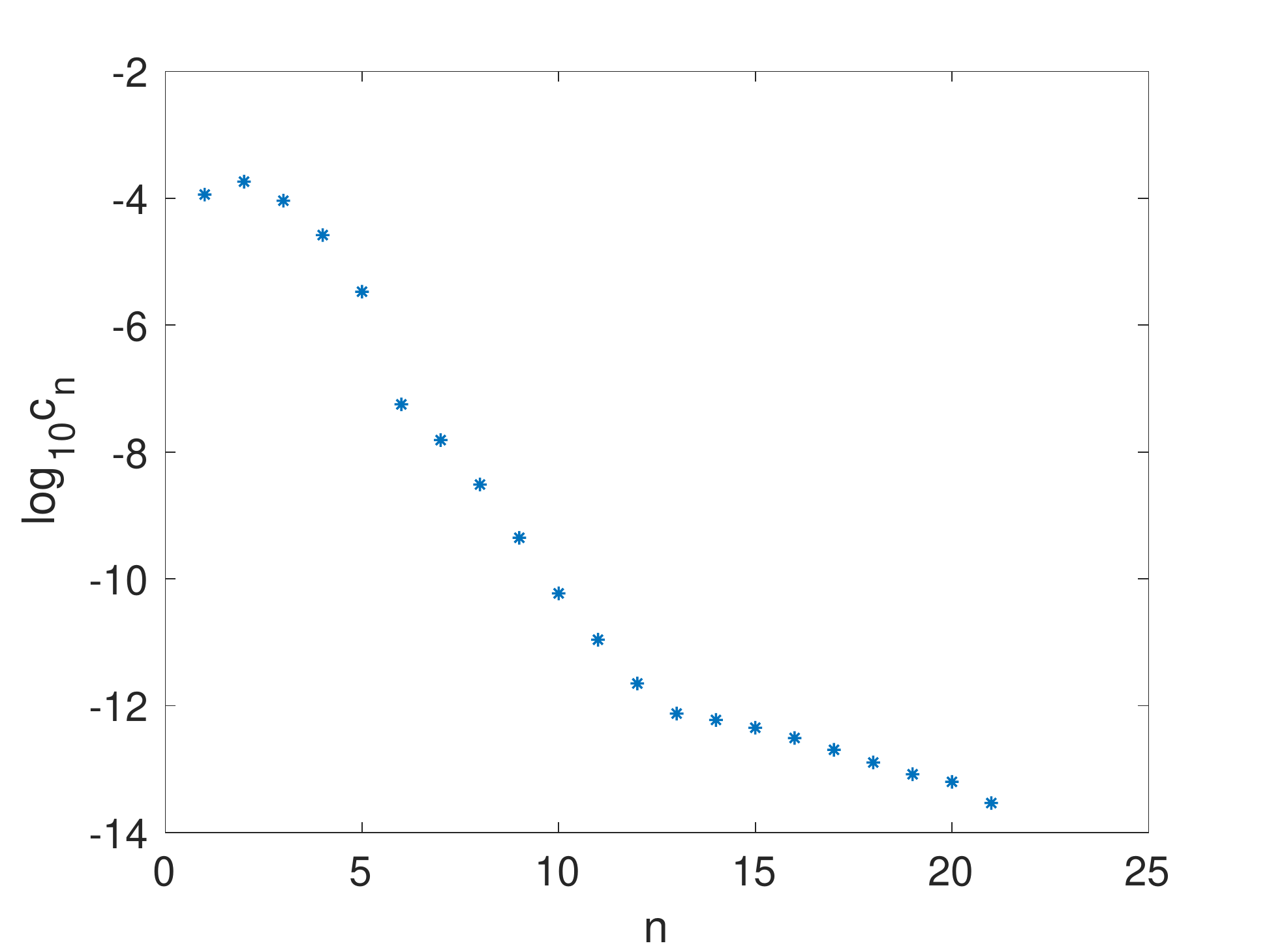}
 \caption{Spectral coefficients $c_{n}$ for the situation in Fig.~\ref{imag} for domain I on the left and for domain~III on the right. }
 \label{coeffimag}
\end{figure}

As already stated, an expansion in terms of Chebyshev polynomials offers the possibility to check the numerical resolution via the decrease of the spectral coefficients $c_{n}$ (\ref{cn}) as computed by an FCT. The coefficients in domain II are shown on the right of Fig.~\ref{sfig}. They decrease to the order of the rounding error for $n\sim 128$. Therefore it would have been possible to use just half of the value of $N_{\rm II}$ without loss of accuracy. The corresponding coefficients for domain~I and III can be seen on the left and right respectively of Fig.~\ref{coeffimag}. Maximal accuracy is reached in this case with roughly 15 polynomials. The spectral coefficients indicate that the solution is as expected smooth on the whole line.

\subsection{Line close to the Stokes' lines}
\begin{figure}[!t]\centering
 \includegraphics[width=0.49\textwidth]{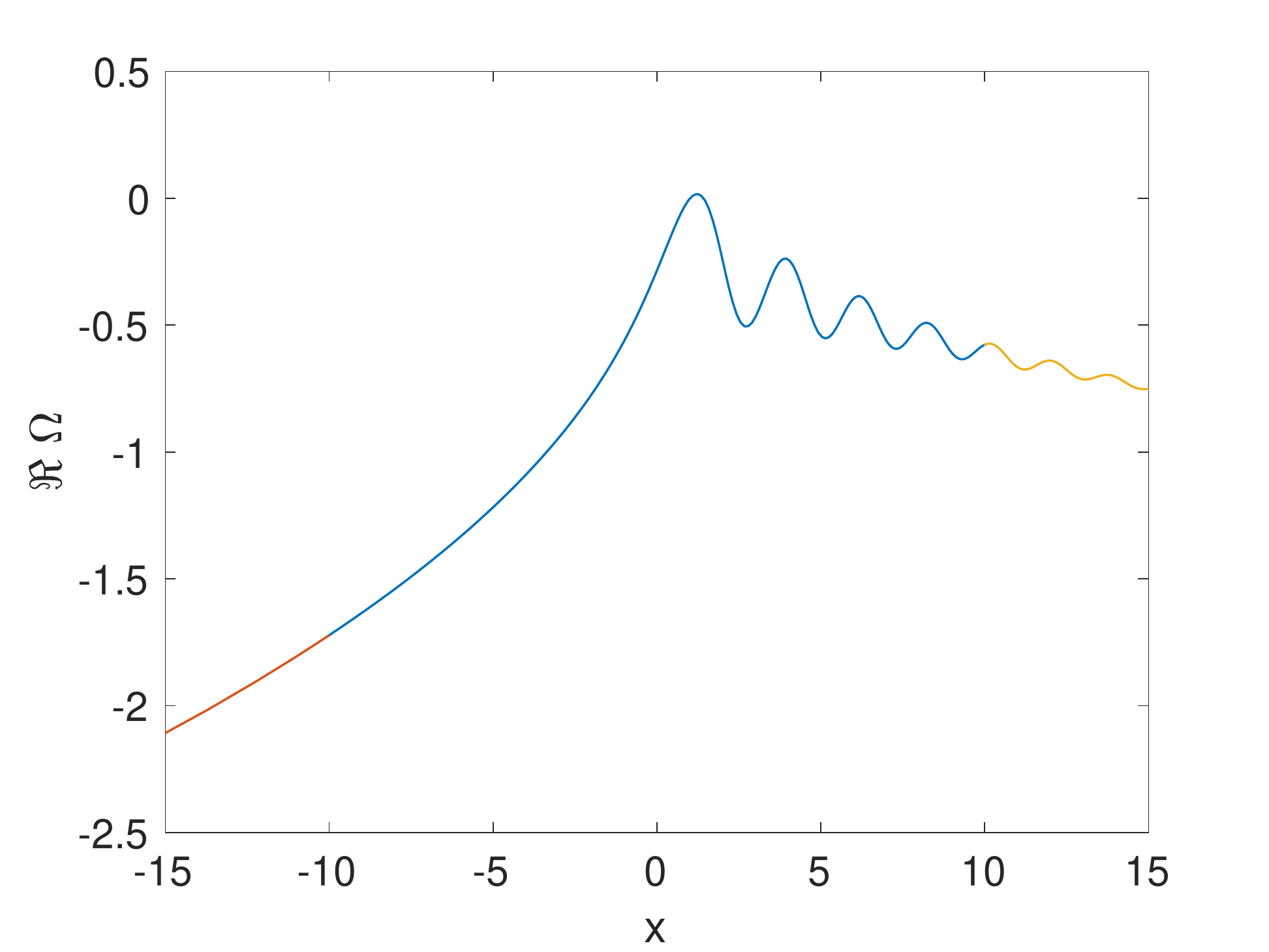}
 \includegraphics[width=0.49\textwidth]{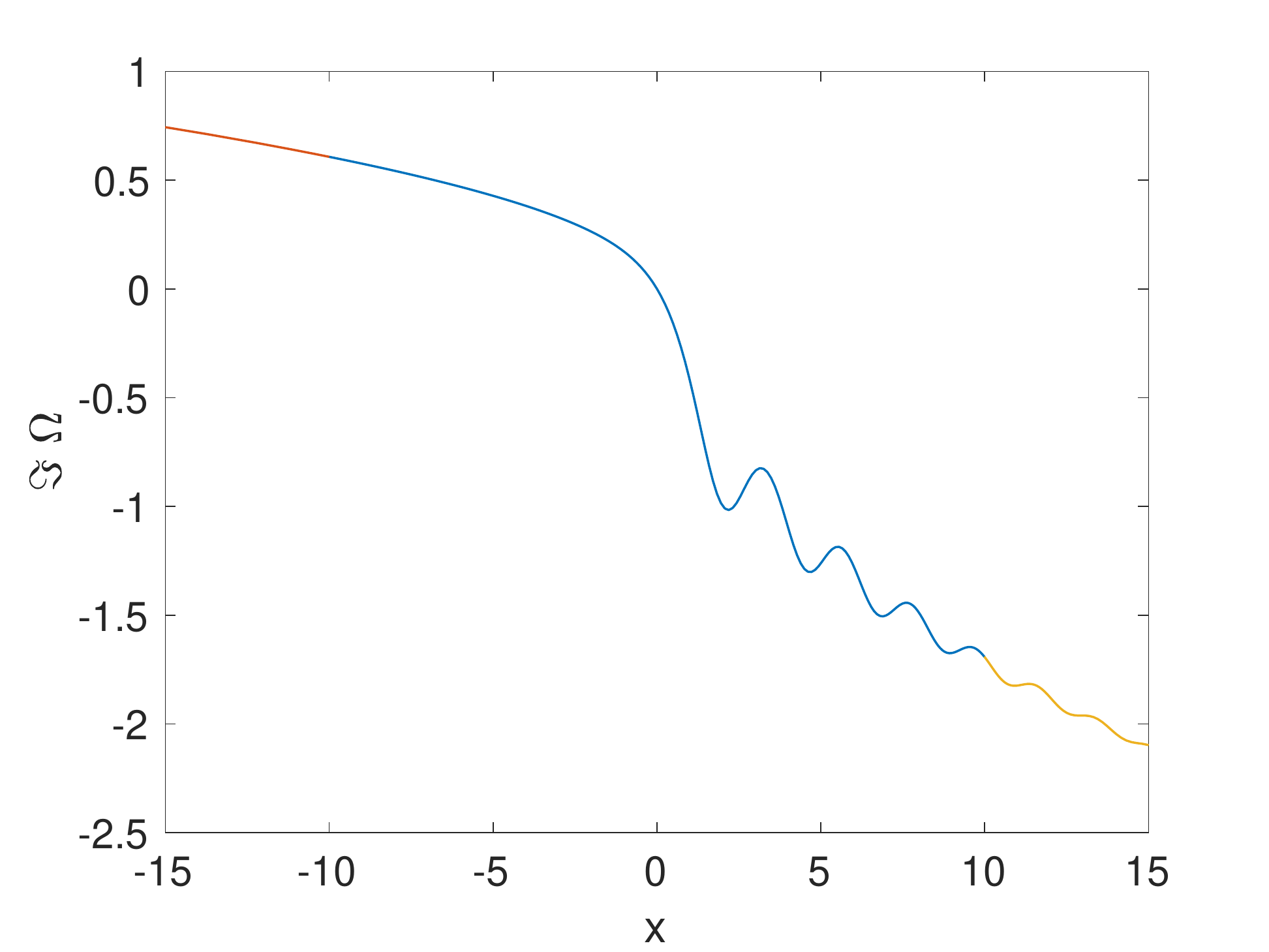}
 \caption{\emph{Tritronqu\'ee} solution to the PI equation on the line $z=\exp(\mathrm{i}(4\pi/5-0.05)x)$ with $x\in\mathbb{R}$, on the left the real part, on the right the imaginary part.} \label{stokes}
\end{figure}

The situation becomes numerically more demanding close to the Stokes' lines $\operatorname{arg}z=\pm 4\pi/5$ because of the oscillatory behavior of the solutions on the latter. We consider the case $z=\exp(\mathrm{i}(4\pi/5-0.05))x$ with $x\in\mathbb{R}$. Again we use $x_{r}=-x_{l}=10$ and $N_{\rm I}=20$, but this time $N_{\rm II}=N_{\rm III}=256$ to provide more resolution near the oscillatory singularity. Though the solution is still obeying the asymptotics (\ref{asym}) on this line, the vicinity to the Stokes' line produces oscillations with decreasing amplitude towards infinity as can be seen in Fig.~\ref{stokes}. Note that there are no oscillations for $x<0$. Again the solutions are smooth on the whole shown interval and thus also at the domain boundaries.

\begin{figure}[!t]\centering
 \includegraphics[width=0.49\textwidth]{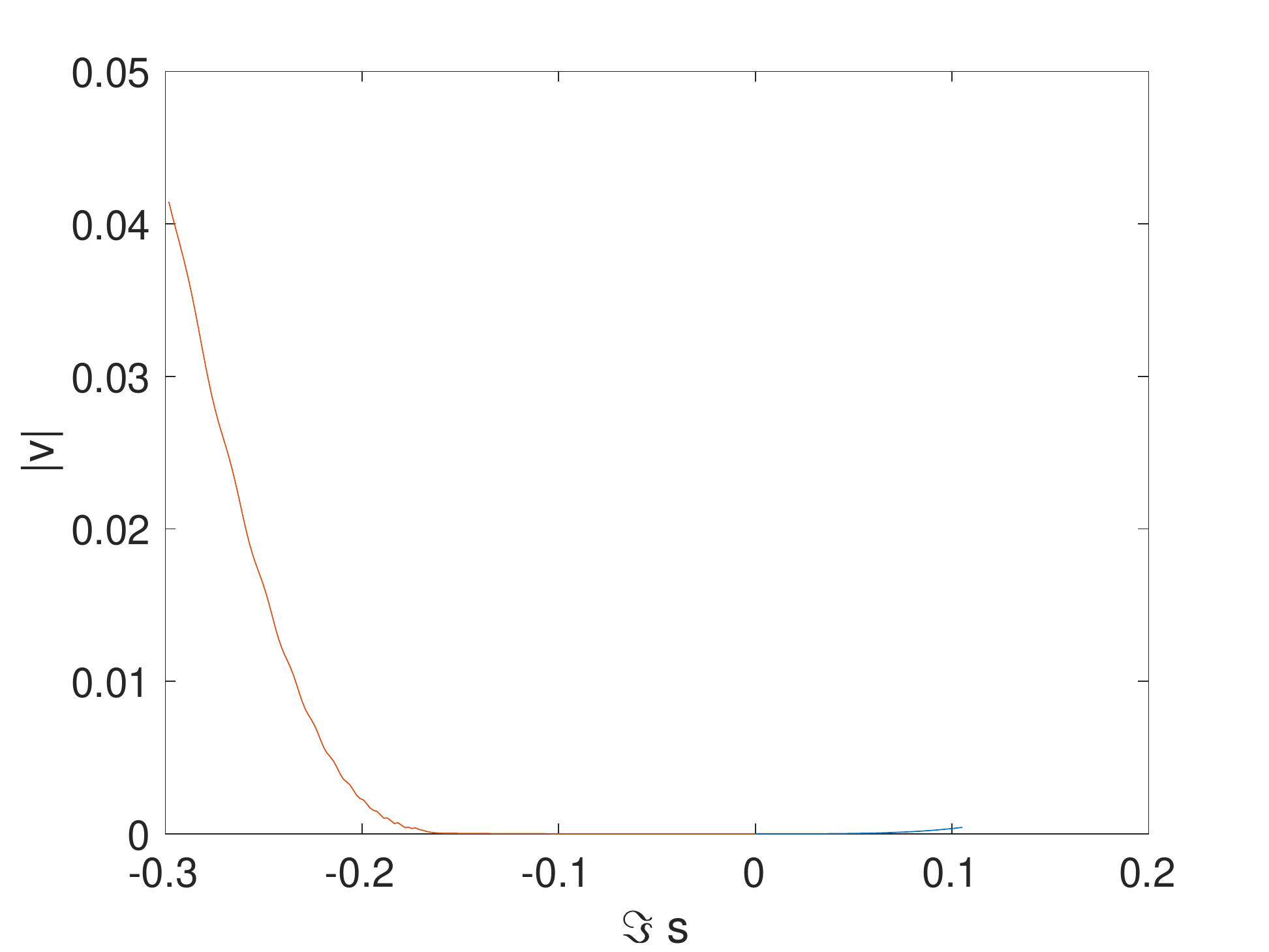}
 \includegraphics[width=0.49\textwidth]{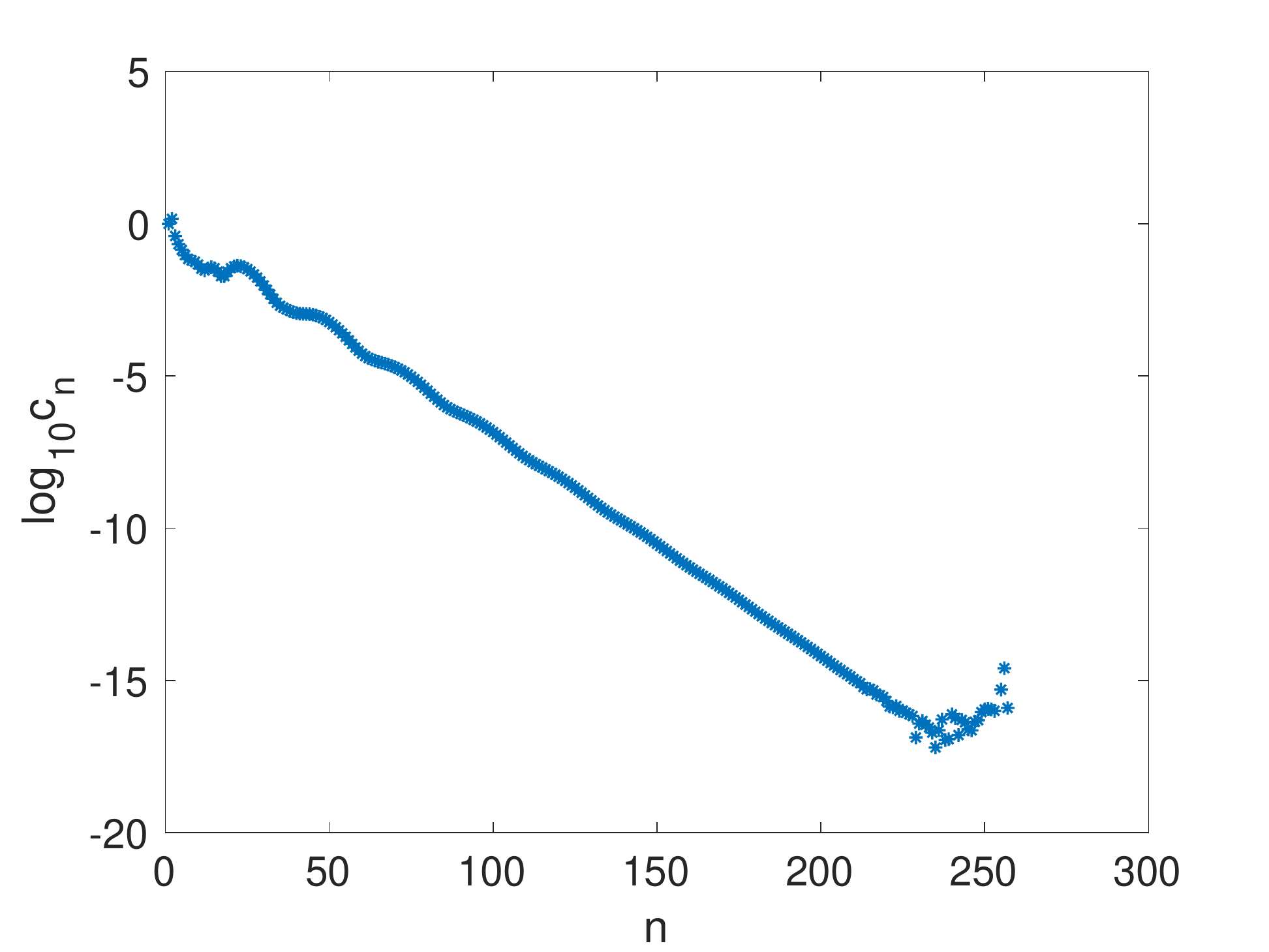}
 \caption{On the left, the computed functions $v$ corresponding to Fig.~\ref{stokes} in domains I (blue) and III (red); on the right the modulus of the spectral coefficients $c_{n}$ in domain II.} \label{stokess}
\end{figure}

The actually computed functions $v$ in domains I and III can be seen in Fig.~\ref{stokess} on the left. In domain I the solution is once more very close to the asymptotic behaviour, the solution is of the order of~$10^{-4}$. However, this is not the case in domain III where the deviation near~$x_{r}$ is of the order~$0.04$.
\begin{figure}[!t]\centering
 \includegraphics[width=0.49\textwidth]{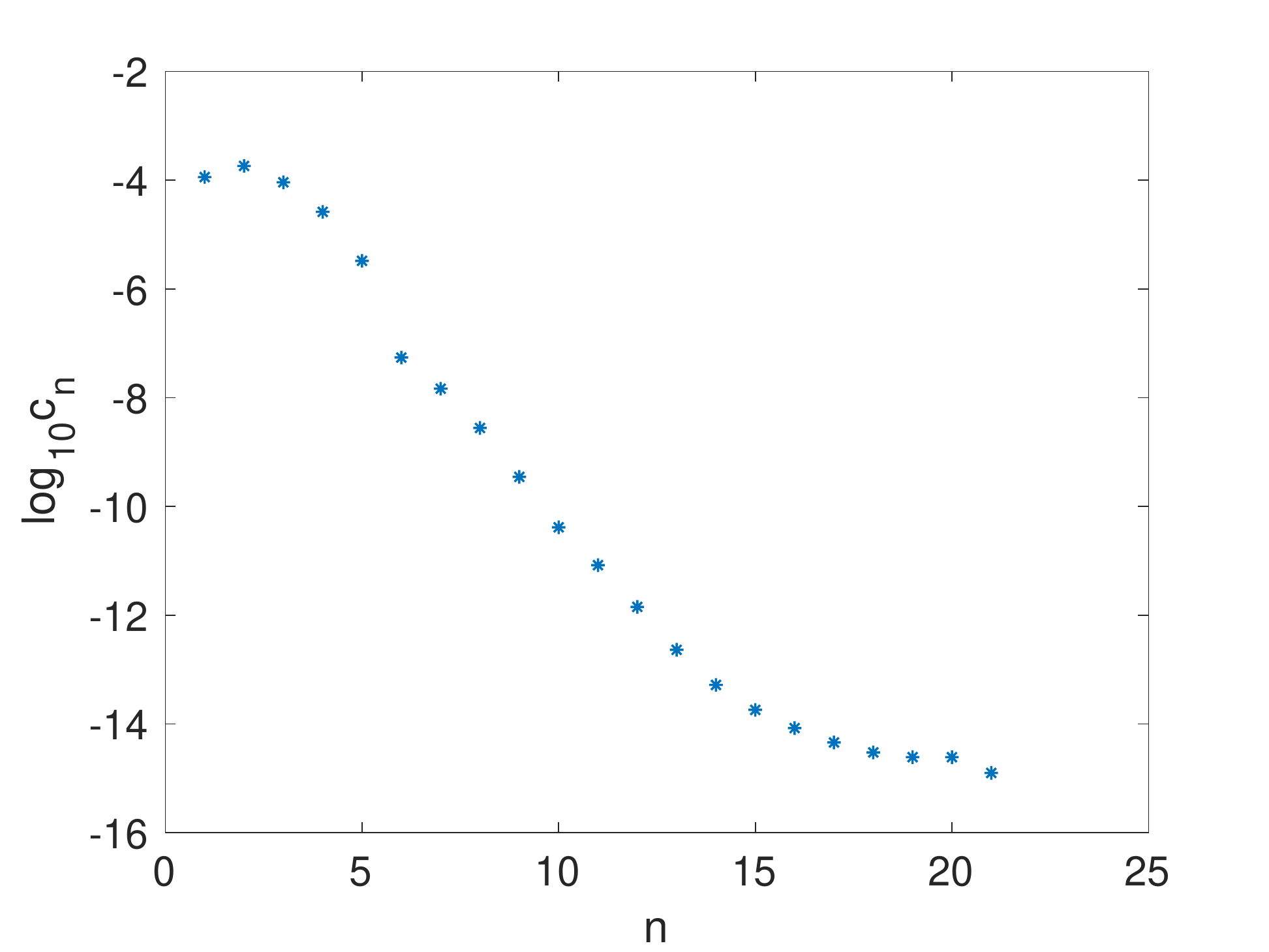}
 \includegraphics[width=0.49\textwidth]{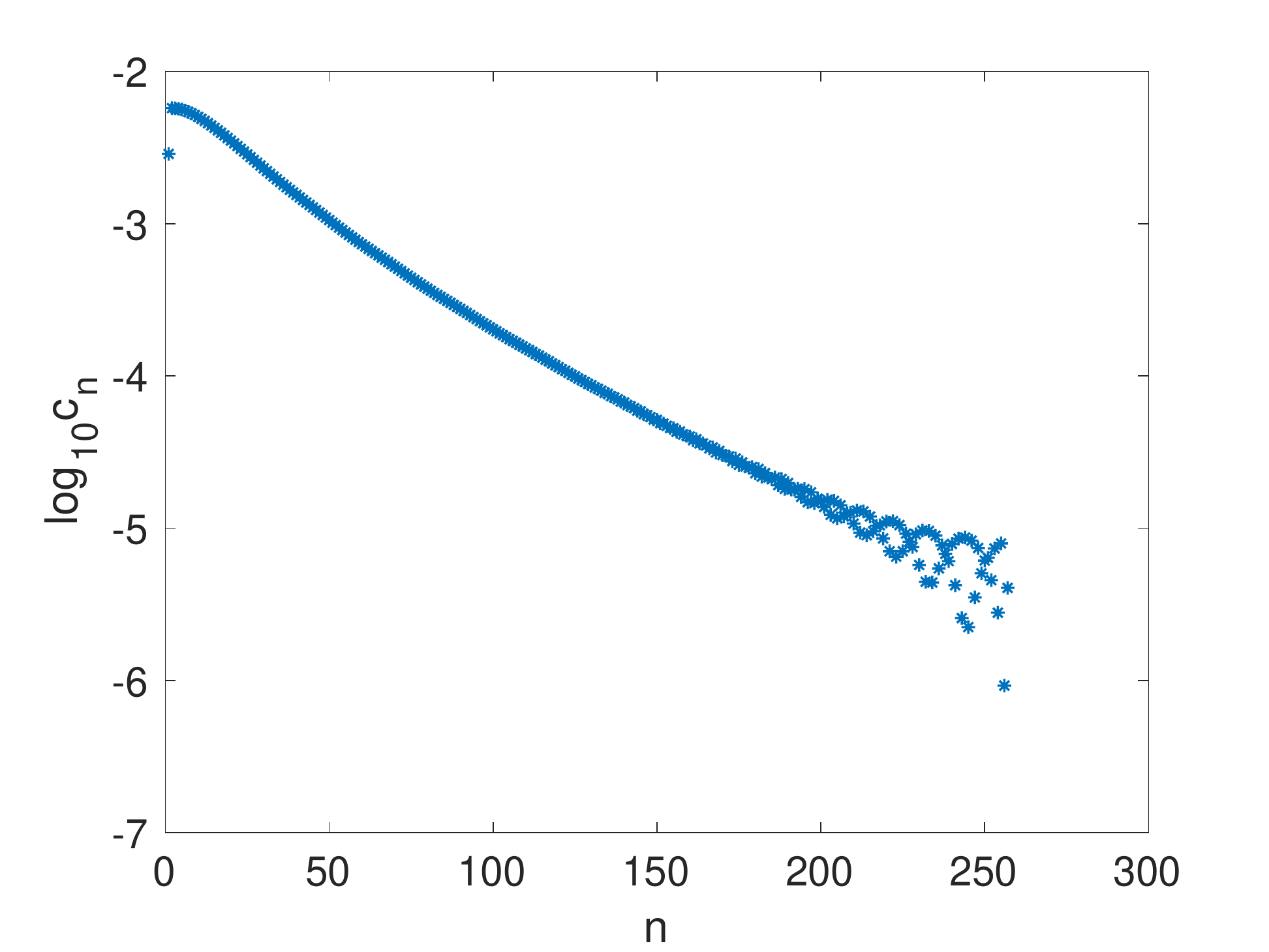}
\caption{Spectral coefficients $c_{n}$ for the situation in Fig.~\ref{stokes} for domain I on the left and for domain III on the right.} \label{stokescoeff}
\end{figure}

The spectral coefficients $c_{n}$ (\ref{cn}) in domain II can be seen on the right of Fig.~\ref{stokess}. Due to the oscillatory nature of the solution, 256 Chebyshev polynomials are needed to reach the level of the rounding error. The spectral coefficients for domains~I and~III are shown in Fig.~\ref{stokescoeff}. As expected the saturation level is reached in domain~I with roughly 15 polynomials. But the vicinity of a~singularity at infinity on the Stokes' line leads to a much smaller decrease of the coefficients in domain III. Though the solution in Fig.~\ref{stokess} appears rather innocent, with 256 polynomials just a~resolution of the order of $10^{-6}$ is reached. Note that this does not change much if a~larger~$x_{r}$ is chosen for the same~$N_{\rm III}$ which increases the numerical resolution. This is a~clear hint on a~loss of regularity in the vicinity of the considered line and possibly even on the line itself near infinity.

\section{Outlook}\label{section4}
In this note we have shown for the example of the \emph{tritronqu\'ee} solutions of PI that Painlev\'e transcendents can be computed on unbounded domains even if they have an algebraic increase in a local parameter near infinity which can be the branch point of a~Riemann surface. A~related example would be the \emph{tritronqu\'ee} solutions of the second equation in the Painlev\'e hierarchy studied in~\cite{GKK}.

If a Painlev\'e transcendent decreases exponentially as the Hastings--McLeod solution~\cite{HM}~$u(x)$ of the Painlev\'e~II equation on the real line for $x\to+\infty$, it is of course sufficient to impose the Dirichlet condition $u(x_{r})=0$ for a~sufficiently large $x_{r}$ such that $u$ vanishes there with numerical accuracy. However, the asymptotic behaviour of the same solution for $x\to-\infty$ is given by $u\sim \sqrt{-x}$. Thus the approach illustrated for the \emph{tritronqu\'ee} solutions could be also applied to the region $x<0$ of the Hastings--McLeod solution.

Since Painlev\'e transcendents can be extended to meromorphic functions in the complex plane, it would be interesting to extend the approach to whole domains of~$\mathbb{C}$. As was pointed out in~\cite{DubrovinGravaKlein}, it is in this case better to set up a~boundary value problem for the Laplace equation where the solution is known to be holomorphic since it is known to be also harmonic there. Since the Laplace equation is linear, no computationally expensive iteration will be needed in this case. In the case of the \emph{tritronqu\'ee} solutions one could think of computing the solution on lines close to the Stokes' lines as in the previous section and use polar coordinates for the Laplace equation to obtain the solution on part of the regular sector. As was mentioned already in~\cite{DubrovinGravaKlein}, the coordinate singularity at the origin leads to a~loss of accuracy there. Such problems can be avoided if as in~\cite{CKSV} for the hypergeometric equation elliptic or rectangular domains are used. It would be also interesting to explore the sector with poles in this way by combining the treatment of the asymptotics of this paper with the use of Pad\'e approximants as in~\cite{FFW, FW}. This will be the subject of further research.

\subsection*{Acknowledgement}
This work was partially supported by the PARI and FEDER programs in 2016 and 2017, by the ANR-FWF project ANuI and by the Marie-Curie RISE network IPaDEGAN. We thank M.~Fasondini for helpful remarks. We are grateful to the anonymous referee for his constructive refereeing and many useful suggestions.

\pdfbookmark[1]{References}{ref}
\LastPageEnding

\end{document}